# The rainbow-spectrum of RNA secondary structures

Thomas J. X. Li · Christian M. Reidys

**Abstract** In this paper we analyze the length-spectrum of rainbows in RNA secondary structures. A rainbow in a secondary structure is a maximal arc with respect to the partial order induced by nesting. We show that there is a significant gap in this length-spectrum. We shall prove that there asymptotically almost surely exists a unique longest rainbow of length at least $n - O(n^{1/2})$ and that with high probability any other rainbow has finite length. We show that the distribution of the length of the longest rainbow converges to a discrete limit law and that, for finite $k$, the distribution of rainbows of length $k$, becomes for large $n$ a negative binomial distribution. We then put the results of this paper into context, comparing the analytical results with those observed in RNA minimum free energy structures, biological RNA structures and relate our findings to the sparsification of folding algorithms.

## 1 Introduction

RNA is a biomolecule involved in a plethora of functions, ranging from catalytic activity to gene expression. A single-stranded RNA molecule has a backbone consisting of nucleotides and can be described by its primary sequence,

Thomas J. X. Li
Biocomplexity Institute of Virginia Tech
Blacksburg, VA 24061, USA
E-mail: thomasli@bi.vt.edu

Christian M. Reidys
Biocomplexity Institute of Virginia Tech
Blacksburg, VA 24061, USA
E-mail: duckcr@bi.vt.edu

i.e., a linear, oriented sequence of the bases $\{\mathbf{A}, \mathbf{U}, \mathbf{G}, \mathbf{C}\}$. In contrast to DNA, an RNA strand folds into a helical configuration of its primary sequence by forming hydrogen bonds between pairs of nucleotides according to Watson-Crick **A-U**, **C-G** and wobble **U-G** base-pairing rules. These structures play a variety of biochemical roles within cells such as: transcription and translation (mRNA links DNA and proteins to convey genetic information with the assistance of tRNA (McCarthy and Holland, 1965)), catalyzing reactions (ribozymes catalyze diverse biological reactions as proteins (Kruger et al, 1982)), gene regulation (miRNA functions in RNA silencing and ncRNA in directing post-transcriptional regulation of gene expression (Eddy, 2001)).

The most prominent class of coarse grained RNA structures are the RNA secondary structures. These are contact structures without any reference of spatial embedding, whose contacts are base pairs subject to certain restrictions. First, their base pairs are canonical pairings: Watson-Crick as well as wobble base pairs. Bonding information such as non-canonical interactions, coaxial stacking of helices, major and minor groove triplexes, and interactions with other molecules are not considered. Secondly, any two base pairs are non-crossing: representing the contact structure as a *diagram*, by drawing its sequence on a horizontal line and each base pair as an arc in the upper half-plane, two arcs $(i_1, j_1)$ and $(i_2, j_2)$ cross if the nucleotides appear in the order $i_1 < i_2 < j_1 < j_2$ in the primary sequence. In this representation, RNA secondary structure contains exclusively non-crossing arcs, see Fig. 1.

The combinatorics of RNA secondary structures was pioneered by Waterman *et al.*, more than three decades ago (Waterman, 1978, 1979; Smith and Waterman, 1978; Howell et al, 1980; Schmitt and Waterman, 1994; Penner and Waterman, 1993). A variety of dynamic programming (DP) algorithms, predicting the minimum free energy (mfe) conformation for RNA molecules, have been derived (Zuker and Sankoff, 1984; Waterman and Smith, 1986; Zuker, 1989; Hofacker et al, 1994). Sparsification is a particular method facilitating a speed up of these DP-routines (Wexler et al, 2007; Salari et al, 2010; Backofen et al, 2011). The method employs the fact that certain matrices of the DP routines are sparse, a fact that greatly simplifies the computation. The theoretical analysis (Wexler et al, 2007) concludes a linear reduction time complexity based on a specific property of arcs in RNA molecules. This property is called *polymer-zeta property* and originates from studies of bonds in proteins. Polymer-zeta asserts that two nucleotides of distance $m$ form a base pair with probability $bm^{-c}$ for some constants $b > 0$, $c > 1$, implying that long-distance base pairs have low probability.

Subsequent analysis revealed that the polymer-zeta property does not hold for general RNA molecules (Backofen et al, 2011), and that sparsification provides only a constant, however significant reduction (Huang and Reidys, 2012).

We shall provide a detailed understanding of the longest, as well as the second-longest arc in RNA secondary structures.

This paper is furthermore motivated by the question of how to interpret the "information" contained in non-coding DNA sequences. In Barrett et al (2017) a sequence-structure correlation of RNA is studied, implying the potential of RNA structure to play a critical role in providing such an interpretation. Accordingly, transcription would not only facilitate the generation of protein (for coding sequences) but also the interpretation of DNA data via forming RNA structure. In other words, it would not be the actual sequence of nucleotides alone but the structures compatible with such sequences that contain crucial information, changing the paradigm of sequence alignments.

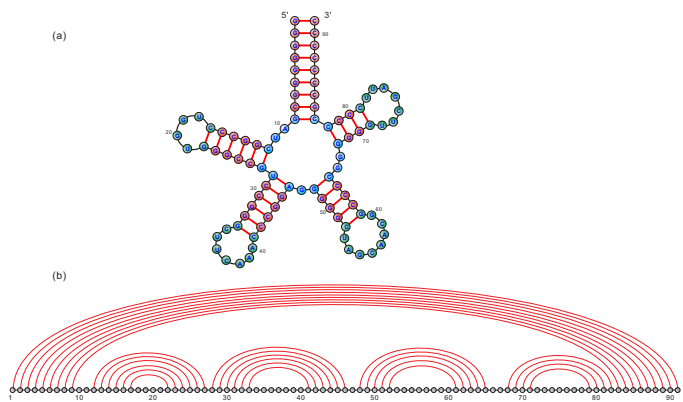

**Fig. 1** An RNA secondary structure represented as a contact graph (a) and as a diagram (b).

In this context it becomes relevant to analyze distances between two paired nucleotides in RNA structures. A particular class of such bonds are rainbows. A rainbow in a secondary structure is a maximal arc with respect to the partial order induced by nesting, i.e. the closing arc of a stem-loop, see Fig. 2. The length of a rainbow $(i,j)$, defined as $j-i$, reflects the size of the corresponding stem-loop. In this paper, we study the length spectrum of rainbows (rainbow-spectrum) in RNA secondary structures.

Rainbows have been studied in Jin and Reidys (2010a,b) in the context of $k$-noncrossing RNA structures. The authors show that the expected number of rainbows is finite and that the endpoint of a rainbow is more likely to occur at the end of the sequence, hinting at the existence of a unique longest rainbow. Another notion closely connected to that of rainbows is the $5'$-$3'$ distance, i.e. the number of rainbows plus the number of unpaired, external nucleotides. The finiteness of the $5'$-$3'$ distance has first been studied in Yoffe et al (2011), where the expected number of rainbows in RNA secondary structures has been obtained. Remarkably, the $5'$-$3'$ distance of biological RNA structures is also observed to be finite, indicating that certain features of random structures can also be observed in biological structures. Han and Reidys (2012) studies rainbows of RNA secondary structures in the context of the $5'$-$3'$ distance. It is shown that this distance satisfies a discrete limit law, implying the finiteness of the $5'$-$3'$ distance of uniformly sampled RNA structures. Clote et al (2012) shows that the expected distance between $5'$ and $3'$ ends of a specific RNA sequence is finite, with respect to the Turner energy model. More importantly, the finiteness of the 5′-3′ distance and the existence of a long rainbow both lead to the effective circularization of linear RNA, which plays an important role in many biological processes (Yoffe et al, 2011).

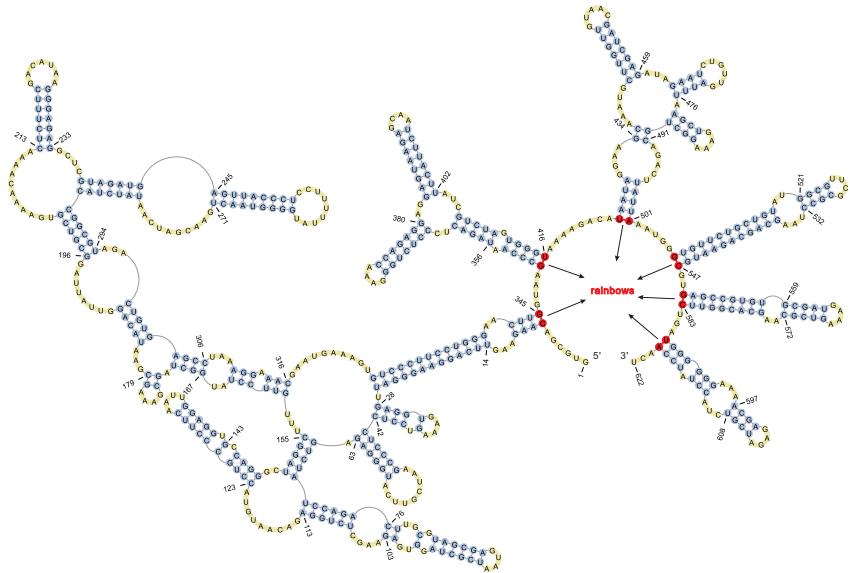

**Fig. 2** The secondary structure of the *P.li.LSUI2* intron (Robart et al, 2014). The structure has six rainbows of lengths 338, 66, 76, 40, 33 and 32, respectively. The figure is generated with the assistance of PseudoViewer3 (Byun and Han, 2009).

The second longest rainbow has in the limit of long sequences, with high probability, finite length. In other words, for any fixed probability, $0 < q < 1$, we find a finite $k(q)$ such that with probability, $q$, a random RNA secondary structure has a second longest rainbow of length at most $k(q)$. However, with probability $o(1)$, there are RNA secondary structures that exhibit a second longest rainbow of order $O(n^{1/2})$ or higher. In fact we shall show that the expected length of the second longest rainbow is $O(n^{1/2})$.

The key results of this paper are the following:

1. in uniformly generated RNA secondary structures the length of the longest rainbow tends, in the limit of long sequences, to a discrete limit law, having an expectation value $n - O(n^{1/2})$. That is, there is a gap in the length-sequence of rainbows, i.e. there exists a unique longest rainbow,
2. with high probability any other rainbow has finite length, $k$,
3. in the limit of long sequences, the distribution of rainbows of length $k$ tends to a negative binomial distribution,
4. mfe-structures also exhibit a unique longest rainbow of order $n - O(n^{1/2})$, and furthermore, with high probability, any other rainbow has finite length.

As for biological structures, in Fig. 2 we display the *P.li.LSUI2* intron (Robart et al, 2014) RNA structure, containing a unique longest rainbow.

In order to obtain the limit distributions of the length of the longest rainbow, we study the generating function of secondary structures having a restricted length of their longest rainbow. This analysis will allow us to compute in Lemma 1 the expectation and variance of the length of the longest rainbow. Having established this we proceed computing the limit distribution in Theorem 3. As for analyzing rainbows of finite length, we consider a bivariate generating function distinguishing the number of rainbows of length $k$ and establish a discrete limit law using the subcritical paradigm (Flajolet and Sedgewick, 2009) of singularity analysis.

This paper is organized as follows: In Section 2, we provide some basic facts of RNA secondary structures. In Section 3, we compute the expectation and variance of the longest rainbow in RNA secondary structures and compute the discrete limit law. In Section 4, we first observe that with high probability we can restrict our analysis to rainbows of finite length and then proceed computing the associated limit distribution. In Section 5, we integrate our results and and discuss them in the context of the $5'$-$3'$ distance and RNA mfe-structures.

## 2 Basic facts

RNA secondary structure can be represented as a *diagram*, a labeled graph over the vertex set $\{1, \ldots, n\}$ whose vertices are arranged in a horizontal line and arcs are drawn in the upper half-plane, see Fig. 1. Clearly, vertices correspond to nucleotides in the primary sequence and arcs correspond to the Watson-Crick as well as wobble base pairs. The *length* of the structure is defined as the number of nucleotides. The length of an arc $(i, j)$ is defined as $j - i$ and an arc of length $k$ is called a $k$-arc. The backbone of a diagram is the sequence of consecutive integers $(1, \ldots, n)$ together with the edges $\{\{i, i+1\} \mid 1 \leq i \leq n-1\}$. We shall distinguish the backbone edge $\{i, i+1\}$ representing a phosphodiester bond, from the arc $(i, i+1)$, which we refer to as a 1-*arc*. Two arcs $(i_1, j_1)$ and $(i_2, j_2)$ are *crossing* if $i_1 < i_2 < j_1 < j_2$. An RNA *secondary structure* is defined as a diagram satisfying the following three conditions (Waterman, 1978):

1. *non-existence of* 1-*arcs*: if $(i, j)$ is an arc, then $j - i \geq 2$,
2. *non-existence of base triples*: any two arcs do not have a common vertex,
3. *non-existence of pseudoknots*: any two arcs are non-crossing, i.e., for two arcs $(i_1, j_1)$ and $(i_2, j_2)$ where $i_1 < i_2$, $i_1 < j_1$, and $i_2 < j_2$, we have either $i_1 < j_1 < i_2 < j_2$ or $i_1 < i_2 < j_2 < j_1$.

A *stack* of length $r$ is a maximal sequence of "parallel" arcs, $((i, j), (i+1, j-1), \ldots, (i+(r-1), j-(r-1)))$. Stacks of length one are energetically unstable and we find typically stacks of length at least two or three in biological structures (Waterman, 1978). A secondary structure, $S$, is $r$-*canonical* if it has minimum stack-length $r$.

Given an RNA secondary structure, $S$, an arc is called a *rainbow* if it is maximal with respect to the partial order $(i, j) \leq (i', j') \iff i' \leq i < j \leq j'$.

I.e. a rainbow is the closing arc of a stem-loop. A secondary structure is called *irreducible* if it contains a rainbow connecting the first and the last vertex in a structure.

We consider RNA secondary structures filtered by *minimum arc-length* and *minimum stack-length*. This filtration is motivated by the fact that for energetic reasons, RNA secondary structures exhibit a minimum arc-length of four and a minimum stack length two or three. The former is a consequence of the rigidity of the molecules backbone (Stein and Waterman, 1979) and the latter a mesomery effect of parallel Watson-Crick or **U-G** base pairs (Hunter and Sanders, 1990; Šponer et al, 2001, 2013).

Let $s_\lambda^{[r]}(n)$ and $f_\lambda^{[r]}(n)$ denote the numbers of $r$-canonical secondary structures and irreducible secondary structures over $n$ nucleotides with minimum arc-length $\lambda$, respectively. We shall simplify notation by writing $s(n)$ and $f(n)$ instead of $s_\lambda^{[r]}(n)$ and $f_\lambda^{[r]}(n)$. The generating functions $\mathbf{S}(x)$ and $\mathbf{F}(x)$ are given by

$$\mathbf{S}(x) = \sum_{n \geq 0} s(n) x^n,$$

$$\mathbf{F}(x) = \sum_{n \geq 1} f(n) x^n,$$

where $f(1) = 1$ represents a single nucleotide, which is irreducible by convention.

These two generating functions have been computed in Waterman (1978); Hofacker et al (1998); Barrett et al (2016).

**Theorem 1** *For any $\lambda, r \in \mathbb{N}$, the generating functions $\mathbf{S}(x)$ and $\mathbf{F}(x)$ satisfy the functional equations*

$$\mathbf{S}(x) = \frac{1}{1 - \mathbf{F}(x)},$$

$$\mathbf{F}(x) - x = \frac{x^{2r}}{1 - x^2} \left( \mathbf{S}(x) - \mathbf{F}(x) + x - \sum_{i=0}^{\lambda-2} x^i \right).$$

*The generating function $\mathbf{S}(x)$ satisfies the functional equation*

$$x^{2r} \mathbf{S}(x)^2 - \mathbf{B}(x) \mathbf{S}(x) + \mathbf{A}(x) = 0,$$

*where*

$$\mathbf{A}(x) = 1 - x^2 + x^{2r},$$

$$\mathbf{B}(x) = (1-x)\mathbf{A}(x) + x^{2r} \sum_{i=0}^{\lambda-2} x^i.$$

*Explicitly, we have*

$$\mathbf{S}(x) = \frac{\mathbf{B}(x) - \sqrt{\mathbf{B}(x)^2 - 4x^{2r} \mathbf{A}(x)}}{2 x^{2r}}.$$

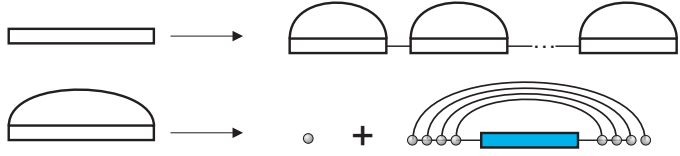

**Fig. 3** The decomposition of a secondary structure and an irreducible structure (reducible structures are colored in blue).

The key idea to prove the functional equation of Theorem 1 is the following: any secondary structures can be decomposed into a sequence of irreducible structures, and any irreducible structure is either a single vertex or the stack containing the rainbow together with the enclosed reducible structure, see Fig. 3.

Singularity analysis of $\mathbf{S}(x)$ (Waterman (1978); Hofacker et al (1998); Barrett et al (2016), implies

**Theorem 2** *For $1 \leq \lambda \leq 4$ and $1 \leq r \leq 3$, the dominant singularity $\rho$ of $\mathbf{F}(x)$ is the minimal positive, real solution of $\mathbf{B}(x)^2 - 4x^{2r}\mathbf{A}(x) = 0$. The singular expansion of $\mathbf{F}(x)$ is given by*

$$\mathbf{F}(x) = \tau + \delta(\rho - x)^{\frac{1}{2}} + \theta(\rho - x) + O\big((\rho - x)^{\frac{3}{2}}\big), \qquad \text{as } x \to \rho,$$

*where $\tau = \mathbf{F}(\rho)$, $\delta$ and $\theta$ are constants, that can be explicitly computed. Furthermore, the coefficients of $\mathbf{F}(x)$ satisfy*

$$[x^n]\mathbf{F}(x) = c\, n^{-\frac{3}{2}} \rho^{-n} \big(1 + O(n^{-1})\big), \qquad \text{as } n \to \infty,$$

*where $c$ is the positive constant $c = -\delta \rho^{\frac{1}{2}} \Gamma(-\frac{1}{2})^{-1}$.*

## 3 The longest rainbow

Our analysis assumes the uniform distribution over all RNA secondary structures of $n$ nucleotides, i.e. the distribution in which each structure has probability $\frac{1}{s(n)}$.

We shall analyze the random variable, $\mathbb{Y}_n$, representing the length of the longest rainbow in an RNA secondary structure of $n$ nucleotides. The generating function of structures, whose rainbows have length less than or equal to $m$ is given by

$$\mathbf{S}_{\leq m+1}(x) = \frac{1}{1 - \mathbf{F}_{\leq m+1}(x)},$$

where $\mathbf{F}_{\leq m}(x) = \sum_{1 \leq i \leq m} f(i) x^i$. By construction, we have

$$\mathbb{P}(\mathbb{Y}_n \leq m) = \frac{[x^n]\mathbf{S}_{\leq m+1}(x)}{[x^n]\mathbf{S}(x)}.$$

In the following we shall derive an asymptotic estimate of $[x^n]\mathbf{S}_{\leq m+1}(x)$. This will imply that the random variable $n - \mathbb{Y}_n$ asymptotically almost surely (a.a.s.) converges to a discrete limit law.

To this end we derive first and second order information about $\mathbb{Y}_n$, which will allow us to apply a large deviation result, instrumental for the proof of our main result.

**Lemma 1** *The expectation and variance of $\mathbb{Y}_n$ are given by*

$$\mathbb{E}[\mathbb{Y}_n] = n - \alpha\, n^{\frac{1}{2}}\big(1 + o(1)\big), \qquad \mathbb{V}[\mathbb{Y}_n] = \beta n^{\frac{3}{2}}\big(1 + o(1)\big), \qquad \text{as } n \to \infty,$$

*where $\alpha = \frac{2\delta \rho^{\frac{1}{2}}}{\sqrt{\pi}(1-\mathbf{F}(\rho))}$ and $\beta = (1 - \frac{\pi}{4})\alpha$ are positive constants.*

*Proof* We consider $\mathbb{P}(\mathbb{Y}_n = n - k)$, by construction, we have

$$\mathbb{P}(\mathbb{Y}_n = n - k) = \frac{[x^n](\mathbf{S}_{\leq n-k+1}(x) - \mathbf{S}_{\leq n-k}(x))}{[x^n]\mathbf{S}(x)}.$$

**Claim 1:** For $k \leq \frac{n}{2}$, we have

$$\mathbb{P}(\mathbb{Y}_n = n - k) = \frac{[x^{k-1}]\Phi'(\mathbf{F}(x))\,[x^{n-k+1}]\mathbf{F}(x)}{[x^n]\mathbf{S}(x)}, \tag{1}$$

where $\Phi(x) = \frac{1}{1-x}$.

*Proof of Claim 1:* The Taylor expansion of $\mathbf{S}_{\leq n-k}(x) = \Phi(\mathbf{F}_{\leq n-k}(x))$ is given by

$$\mathbf{S}_{\leq n-k}(x) = \Phi(\mathbf{F}_{\leq n-k}(x)) = \sum_{i \geq 0} \frac{\Phi^{(i)}(\mathbf{F}(x))}{i!}(\mathbf{F}_{\leq n-k}(x) - \mathbf{F}(x))^i. \tag{2}$$

Note that $[x^n](\mathbf{F}_{\leq n-k}(x) - \mathbf{F}(x))^i = 0$ for $i \geq 2$, since $k \leq \frac{n}{2}$ and $\deg(\mathbf{F}_{\leq n-k}(x) - \mathbf{F}(x)) > \frac{n}{2}$. By taking the coefficient of $x^n$ in eq. (2), we obtain

$$[x^n]\mathbf{S}_{\leq n-k}(x) = [x^n]\Big(\Phi(\mathbf{F}(x)) + \Phi'(\mathbf{F}(x))(\mathbf{F}_{\leq n-k}(x) - \mathbf{F}(x))\Big). \tag{3}$$

Similarly, eq. (3) holds for $[x^n]\mathbf{S}_{\leq n-k+1}(x)$. Therefore, we arrive at

$$\begin{aligned}\mathbb{P}(\mathbb{Y}_n = n - k) &= \frac{[x^n](\Phi'(\mathbf{F}(x))(\mathbf{F}_{\leq n-k+1}(x) - \mathbf{F}_{\leq n-k}(x)))}{[x^n]\mathbf{S}(x)} \\ &= \frac{[x^n](\Phi'(\mathbf{F}(x))f(n-k+1)x^{n-k+1})}{[x^n]\mathbf{S}(x)} \\ &= \frac{[x^{k-1}]\Phi'(\mathbf{F}(x))\,[x^{n-k+1}]\mathbf{F}(x)}{[x^n]\mathbf{S}(x)}.\end{aligned}$$

∎

The rainbow-spectrum of RNA secondary structures9

**Claim 2:**

$$\sum_{1\leq k\leq \frac{n}{2}} (k-1)\mathbb{P}(\mathbb{Y}_n = n-k) = \alpha\, n^{\frac{1}{2}}\bigl(1+o(1)\bigr), \qquad \text{as } n\to\infty. \qquad (4)$$

*Proof of Claim 2:* We shall first derive an estimate of $\mathbb{P}(\mathbb{Y}_n = n-k)$ from Claim 1. By Theorem 2, we have the singular expansions of $\mathbf{F}(x)$, $\mathbf{S}(x) = \Phi(\mathbf{F}(x))$ and $\Phi'(\mathbf{F}(x))$

$$\mathbf{F}(x) = \tau + \delta(\rho-x)^{\frac{1}{2}} + \theta(\rho-x) + O\bigl((\rho-x)^{\frac{3}{2}}\bigr),$$

$$\Phi(\mathbf{F}(x)) = \Phi(\tau) + \Phi'(\tau)\delta(\rho-x)^{\frac{1}{2}} + \theta_1(\rho-x) + O\bigl((\rho-x)^{\frac{3}{2}}\bigr),$$

$$\Phi'(\mathbf{F}(x)) = \Phi'(\tau) + \Phi''(\tau)\delta(\rho-x)^{\frac{1}{2}} + \theta_2(\rho-x) + O\bigl((\rho-x)^{\frac{3}{2}}\bigr),$$

where $\theta, \theta_1, \theta_2$ are constants, and the singular expansions of $\Phi(\mathbf{F}(x))$ and $\Phi'(\mathbf{F}(x))$ are obtained by combining the regular expansions of $\Phi(x)$ and $\Phi'(x)$ with the singular expansion of $\mathbf{F}(x)$ (the subcritical case, see Flajolet and Sedgewick (2009) pp. 411). The Transfer Theorem (Flajolet and Sedgewick (2009) pp. 390) then implies

$$\begin{aligned}
\frac{[x^{n-k+1}]\mathbf{F}(x)}{[x^n]\mathbf{S}(x)} &= \frac{(n-k+1)^{-\frac{3}{2}}\rho^{-n+k-1}\bigl(1+O((n-k)^{-1})\bigr)}{\Phi'(\tau)n^{-\frac{3}{2}}\rho^{-n}\bigl(1+O(n^{-1})\bigr)} \\
&= \frac{(1-\frac{k-1}{n})^{-\frac{3}{2}}\rho^{k-1}\bigl(1+O(n^{-1})\bigr)}{\Phi'(\tau)\bigl(1+O(n^{-1})\bigr)} \qquad (5)\\
&= \frac{(1-\frac{k-1}{n})^{-\frac{3}{2}}\rho^{k-1}}{\Phi'(\tau)}\bigl(1+O(n^{-1})\bigr), \quad \text{as } n\to\infty,\ k\leq \frac{n}{2},
\end{aligned}$$

$$[x^{k-1}]\Phi'(\mathbf{F}(x)) = \frac{\delta\rho^{\frac{1}{2}}\Phi''(\tau)}{-\Gamma(-\frac{1}{2})}(k-1)^{-\frac{3}{2}}\rho^{-k+1}\bigl(1+O(k^{-1})\bigr), \quad \text{as } k\to\infty. \quad (6)$$

Inserting this into eq. (1) and using $\tau = \mathbf{F}(\rho)$, we obtain

$$\begin{aligned}
&\mathbb{P}(\mathbb{Y}_n = n-k) \\
&= [x^{k-1}]\Phi'(\mathbf{F}(x))\frac{[x^{n-k+1}]\mathbf{F}(x)}{[x^n]\mathbf{S}(x)} \qquad (7)\\
&= -\frac{\delta\rho^{\frac{1}{2}}\Phi''(\tau)}{\Gamma(-\frac{1}{2})\Phi'(\tau)}\left(1-\frac{k-1}{n}\right)^{-\frac{3}{2}}(k-1)^{-\frac{3}{2}}\bigl(1+O(k^{-1})\bigr)\bigl(1+O(n^{-1})\bigr),
\end{aligned}$$

as $k\to\infty$, $n\to\infty$ and $k\leq \frac{n}{2}$. In view of the fact that the probability $\mathbb{P}(\mathbb{Y}_n = n-k)$ is at most 1, we have $\sum_{1\leq k\leq n^{\frac{1}{8}}}(k-1)\mathbb{P}(\mathbb{Y}_n = n-k) = O(n^{\frac{1}{8}}\cdot n^{\frac{1}{8}}) = o(n^{\frac{1}{2}})$. Furthermore for large $k$, we have eq. (7). This motivates to split the summation of eq. (4) and to consider the term $\sum_{n^{\frac{1}{8}} \leq k \leq \frac{n}{2}} (k-1)\mathbb{P}(\mathbb{Y}_n = n-k)$ separately, as this allows to employ eq. (7). This leads to

$$
\begin{aligned}
\sum_{n^{\frac{1}{8}} \leq k \leq \frac{n}{2}} & (k-1)\mathbb{P}(\mathbb{Y}_n = n-k) \\
&= \frac{\delta \rho^{\frac{1}{2}} \Phi''(\tau)}{-\Gamma(-\frac{1}{2})\Phi'(\tau)} \sum_{n^{\frac{1}{8}} \leq k \leq \frac{n}{2}} \left(1 - \frac{k-1}{n}\right)^{-\frac{3}{2}} (k-1)^{-\frac{1}{2}} \left(1 + O(k^{-1})\right)\left(1 + O(n^{-1})\right) \\
&= \frac{\delta \rho^{\frac{1}{2}} \Phi''(\tau)}{-\Gamma(-\frac{1}{2})\Phi'(\tau)} n^{\frac{1}{2}} \sum_{n^{\frac{1}{8}} \leq k \leq \frac{n}{2}} \left(1 - \frac{k-1}{n}\right)^{-\frac{3}{2}} \left(\frac{k-1}{n}\right)^{-\frac{1}{2}} \left(\frac{1}{n}\right) \left(1 + O(k^{-1})\right)\left(1 + O(n^{-1})\right) \\
&= \frac{\delta \rho^{\frac{1}{2}} \Phi''(\tau)}{-\Gamma(-\frac{1}{2})\Phi'(\tau)} n^{\frac{1}{2}} \int_0^{\frac{1}{2}} (1-x)^{-\frac{3}{2}} x^{-\frac{1}{2}} \mathrm{d}x \left(1 + o(1)\right)\left(1 + O(n^{-1})\right) \\
&= \alpha n^{\frac{1}{2}} \left(1 + o(1)\right), \qquad \text{as } n \to \infty,
\end{aligned}
\tag{8}
$$

where $\alpha$ is given by

$$
\alpha = \frac{\delta \rho^{\frac{1}{2}} \Phi''(\tau)}{-\Gamma(-\frac{1}{2})\Phi'(\tau)} \int_0^{\frac{1}{2}} (1-x)^{-\frac{3}{2}} x^{-\frac{1}{2}} \mathrm{d}x = \frac{2\delta \rho^{\frac{1}{2}}}{\sqrt{\pi}(1 - \mathbf{F}(\rho))}.
$$

To see the third equality in eq. (8), we first derive

$$
\begin{aligned}
\sum_{n^{\frac{1}{8}} \leq k \leq \frac{n}{2}} & \left(1 - \frac{k-1}{n}\right)^{-\frac{3}{2}} \left(\frac{k-1}{n}\right)^{-\frac{1}{2}} \left(\frac{1}{n}\right) \\
&= \sum_{1 \leq k \leq \frac{n}{2}} \left(1 - \frac{k-1}{n}\right)^{-\frac{3}{2}} \left(\frac{k-1}{n}\right)^{-\frac{1}{2}} \left(\frac{1}{n}\right) - \sum_{1 \leq k \leq n^{\frac{1}{8}}} \left(1 - \frac{k-1}{n}\right)^{-\frac{3}{2}} \left(\frac{k-1}{n}\right)^{-\frac{1}{2}} \left(\frac{1}{n}\right) \\
&= \int_0^{\frac{1}{2}} (1-x)^{-\frac{3}{2}} x^{-\frac{1}{2}} \mathrm{d}x \left(1 + o(1)\right) - \int_0^{\frac{n^{\frac{1}{8}} - 1}{n}} (1-x)^{-\frac{3}{2}} x^{-\frac{1}{2}} \mathrm{d}x \left(1 + o(1)\right) \qquad (9) \\
&= \int_0^{\frac{1}{2}} (1-x)^{-\frac{3}{2}} x^{-\frac{1}{2}} \mathrm{d}x \left(1 + o(1)\right) + O(n^{-\frac{7}{16}}) \\
&= \int_0^{\frac{1}{2}} (1-x)^{-\frac{3}{2}} x^{-\frac{1}{2}} \mathrm{d}x \left(1 + o(1)\right), \qquad \text{as } n \to \infty.
\end{aligned}
$$

Here we can estimate the two sums by the integrals in eq. (9), since the integrals converge. The error term $\sum_{n^{\frac{1}{8}} \leq k \leq \frac{n}{2}} \left(1 - \frac{k-1}{n}\right)^{-\frac{3}{2}} \left(\frac{k-1}{n}\right)^{-\frac{1}{2}} \left(\frac{1}{n}\right) O(k^{-1}) = O\left(\sum_{n^{\frac{1}{8}} \leq k \leq \frac{n}{2}} \left(1 - \frac{k-1}{n}\right)^{-\frac{3}{2}} \left(\frac{k-1}{n}\right)^{-\frac{3}{2}} \left(\frac{1}{n}\right)^2\right)$, as $n \to \infty$. Furthermore, we can show that $\sum_{n^{\frac{1}{8}} \leq k \leq \frac{n}{2}} \left(1 - \frac{k-1}{n}\right)^{-\frac{3}{2}} \left(\frac{k-1}{n}\right)^{-\frac{3}{2}} \left(\frac{1}{n}\right)^2 = O(n^{-\frac{9}{16}})$. As a result we

obtain

$$\sum_{n^{\frac{1}{8}} \leq k \leq \frac{n}{2}} \left(1 - \frac{k-1}{n}\right)^{-\frac{3}{2}} \left(\frac{k-1}{n}\right)^{-\frac{1}{2}} \left(\frac{1}{n}\right) (1 + O(k^{-1})) = \int_0^{\frac{1}{2}} (1-x)^{-\frac{3}{2}} x^{-\frac{1}{2}} dx \, (1 + o(1)),$$

as $n \to \infty$.

Combining eq. (8) and $\sum_{1 \leq k \leq n^{\frac{1}{8}}} (k-1) \mathbb{P}(\mathbb{Y}_n = n - k) = o(n^{\frac{1}{2}})$, we derive eq. (4). ∎

**Claim 3:**

$$\sum_{\frac{n}{2} < k \leq n} (k-1) \mathbb{P}(\mathbb{Y}_n = n - k) = o(n^{\frac{1}{2}}), \qquad \text{as } n \to \infty.$$

*Proof of Claim 3:* We compute

$$\sum_{\frac{n}{2} < k \leq n} (k-1) \mathbb{P}(\mathbb{Y}_n = n - k)$$

$$\leq n \sum_{\frac{n}{2} < k \leq n} \mathbb{P}(\mathbb{Y}_n = n - k)$$

$$= n \left(1 - \sum_{1 \leq k \leq \frac{n}{2}} \mathbb{P}(\mathbb{Y}_n = n - k)\right)$$

$$= n \left(1 - \sum_{1 \leq k < n^{\frac{2}{5}}} \mathbb{P}(\mathbb{Y}_n = n - k) - \sum_{n^{\frac{2}{5}} \leq k \leq \frac{n}{2}} \mathbb{P}(\mathbb{Y}_n = n - k)\right). \qquad (10)$$

We choose $k = n^{\frac{2}{5}}$ as the cutoff, in order to employ eq. (7) for large $k$ and as a result the error term of the estimate is of order $o(n^{\frac{1}{2}})$. We proceed by computing $\sum_{1 \leq k < n^{\frac{2}{5}}} \mathbb{P}(\mathbb{Y}_n = n - k)$ and $\sum_{n^{\frac{2}{5}} \leq k \leq \frac{n}{2}} \mathbb{P}(\mathbb{Y}_n = n - k)$.

For any $1 \leq k < n^{\frac{2}{5}}$, we derive from eq. (1) together with eq. (5)

$$\mathbb{P}(\mathbb{Y}_n = n - k) = [x^{k-1}] \Phi'(\mathbf{F}(x)) \frac{[x^{n-k+1}] \mathbf{F}(x)}{[x^n] \mathbf{S}(x)}$$

$$= \frac{\rho^{k-1} [x^{k-1}] \Phi'(\mathbf{F}(x))}{\Phi'(\tau)} \left(1 - \frac{k-1}{n}\right)^{-\frac{3}{2}} (1 + O(n^{-1})) \qquad (11)$$

$$= c \, b_k \, \rho^{k-1} \left(1 + O(n^{-\frac{3}{5}})\right) \left(1 + O(n^{-1})\right)$$

$$= c \, b_k \, \rho^{k-1} (1 + o(n^{-\frac{1}{2}})), \qquad \text{as } n \to \infty,$$

where $c = \Phi'(\mathbf{F}(\rho))^{-1}$, $b_k = [x^{k-1}] \Phi'(\mathbf{F}(x))$. The third equation follows from $\lim_{n \to \infty} \left(1 - \frac{k-1}{n}\right)^{-\frac{3}{2}} = 1 + O(n^{-\frac{3}{5}})$, since $k < n^{\frac{2}{5}}$. Thus, inserting eq. (11)

into the first sum in eq. (10), we obtain

$$\sum_{1 \leq k < n^{\frac{2}{5}}} \mathbb{P}(\mathbb{Y}_n = n - k)$$

$$= \sum_{1 \leq k < n^{\frac{2}{5}}} c\, b_k\, \rho^{k-1} \bigl(1 + o(n^{-\frac{1}{2}})\bigr)$$

$$= \Bigl(c \sum_{k \geq 1} b_k\, \rho^{k-1} - c \sum_{k \geq n^{\frac{2}{5}}} b_k\, \rho^{k-1}\Bigr)\bigl(1 + o(n^{-\frac{1}{2}})\bigr)$$

$$= \Bigl(1 - c \sum_{k \geq n^{\frac{2}{5}}} b_k\, \rho^{k-1}\Bigr)\bigl(1 + o(n^{-\frac{1}{2}})\bigr) \tag{12}$$

$$= \Bigl(1 - \frac{c\,\delta\,\rho^{\frac{1}{2}}\Phi''(\tau)}{-\Gamma(-\frac{1}{2})} \sum_{k \geq n^{\frac{2}{5}}} k^{-\frac{3}{2}}\bigl(1 + O(k^{-1})\bigr)\Bigr)\bigl(1 + o(n^{-\frac{1}{2}})\bigr) \tag{13}$$

$$= \Bigl(1 - \frac{2c\,\delta\,\rho^{\frac{1}{2}}\Phi''(\tau)}{-\Gamma(-\frac{1}{2})}\, n^{-\frac{1}{5}}\bigl(1 + O(n^{-\frac{2}{5}})\bigr)\Bigr)\bigl(1 + o(n^{-\frac{1}{2}})\bigr) \tag{14}$$

$$= \Bigl(1 - \alpha n^{-\frac{1}{5}}\bigl(1 + O(n^{-\frac{2}{5}})\bigr)\Bigr)\bigl(1 + o(n^{-\frac{1}{2}})\bigr) \tag{15}$$

$$= 1 - \alpha n^{-\frac{1}{5}} + o(n^{-\frac{1}{2}}), \qquad \text{as } n \to \infty.$$

Eq. (12) follows from $\sum_{k \geq 1} b_k\, \rho^{k-1} = \Phi'(\mathbf{F}(\rho)) = c^{-1}$ and eq. (13) follows from eq. (6), since $k \geq n^{\frac{2}{5}}$ tends to infinity. For eq. (14), we know $\sum_{k \geq n^{\frac{2}{5}}} k^{-\frac{3}{2}} = \zeta(3/2, n^{\frac{2}{5}})$, where $\zeta(s,n) = \sum_{i=0}^{\infty} (n+i)^{-s}$ is the Hurwitz-Zeta function. It is well known that, for $s > 1$, as real number $n \to \infty$, the Hurwitz-Zeta function has the asymptotic expansion $\zeta(s,n) = \frac{n^{1-s}}{s-1}\bigl(1 + O(n^{-1})\bigr)$. Then we derive $\sum_{k \geq n^{\frac{2}{5}}} k^{-\frac{3}{2}} = \zeta(3/2, n^{\frac{2}{5}}) = 2n^{-\frac{1}{5}}\bigl(1 + O(n^{-\frac{2}{5}})\bigr)$. Similarly, the error term $\sum_{k \geq n^{\frac{2}{5}}} k^{-\frac{3}{2}} O(k^{-1}) = O(\zeta(5/2, n^{\frac{2}{5}})) = O(n^{-\frac{3}{5}})$. Thus we derive $\sum_{k \geq n^{\frac{2}{5}}} k^{-\frac{3}{2}}\bigl(1 + O(k^{-1})\bigr) = 2n^{-\frac{1}{5}}\bigl(1 + O(n^{-\frac{2}{5}})\bigr)$. Eq. (15) follows from the definition of $\alpha = -\frac{2c\,\delta\,\rho^{\frac{1}{2}}\Phi''(\tau)}{\Gamma(-\frac{1}{2})}$.

As for the second sum, we have

$$\sum_{n^{\frac{2}{5}} \leq k \leq \frac{n}{2}} \mathbb{P}(\mathbb{Y}_n = n - k)$$

$$= \frac{\delta \rho^{\frac{1}{2}} \Phi''(\tau)}{-\Gamma(-\frac{1}{2})\Phi'(\tau)} \sum_{n^{\frac{2}{5}} \leq k \leq \frac{n}{2}} \left(1 - \frac{k-1}{n}\right)^{-\frac{3}{2}} (k-1)^{-\frac{3}{2}} \left(1 + O(k^{-1})\right)\left(1 + O(n^{-1})\right)$$
(16)

$$= \frac{\alpha}{2} \sum_{n^{\frac{2}{5}} \leq k \leq \frac{n}{2}} \left(1 - \frac{k-1}{n}\right)^{-\frac{3}{2}} (k-1)^{-\frac{3}{2}} \left(1 + O(k^{-1})\right)\left(1 + O(n^{-1})\right)$$

$$= \frac{\alpha}{2} \cdot 2n^{-\frac{1}{5}} (1 + O(n^{-\frac{2}{5}})) \left(1 + O(n^{-1})\right)$$
(17)

$$= \alpha n^{-\frac{1}{5}} + o(n^{-\frac{1}{2}}), \qquad \text{as } n \to \infty.$$

Eq. (16) follows from eq. (7), as $k \to \infty$, $n \to \infty$ and $k \leq \frac{n}{2}$. In eq. (17), the summation is approximated by the Euler-Maclaurin summation formula (see, for example, Graham and Knuth and Patashnik (1994)) as follows

$$\sum_{n^{\frac{2}{5}} \leq k \leq \frac{n}{2}} \left(1 - \frac{k-1}{n}\right)^{-\frac{3}{2}} (k-1)^{-\frac{3}{2}}$$

$$= \int_{n^{\frac{2}{5}}}^{\frac{n}{2}} f(x) \mathrm{d}x - \frac{1}{2} f(x) \Big|_{n^{\frac{2}{5}}}^{\frac{n}{2}} + \frac{1}{12} f'(x) \Big|_{n^{\frac{2}{5}}}^{\frac{n}{2}} + R_2(n)$$

$$= 2n^{-\frac{1}{5}}(1 + O(n^{-\frac{3}{5}})) - \frac{1}{2} O(n^{-\frac{3}{5}}) + \frac{1}{12} O(n^{-1}) + O(n^{-1})$$

$$= 2n^{-\frac{1}{5}}(1 + O(n^{-\frac{2}{5}})), \qquad \text{as } n \to \infty,$$

where $f(x) = (1 - \frac{x}{n})^{-\frac{3}{2}} x^{-\frac{3}{2}}$ and the remainder $R_2(n) = O\left((2\pi)^{-2} \int_{n^{\frac{2}{5}}}^{\frac{n}{2}} |f''(x)| \mathrm{d}x\right) = O(n^{-1})$. Similarly, the error term satisfies $\sum_{n^{\frac{2}{5}} \leq k \leq \frac{n}{2}} \left(1 - \frac{k-1}{n}\right)^{-\frac{3}{2}} (k-1)^{-\frac{3}{2}} O(k^{-1}) = O(n^{-\frac{3}{5}})$. Accordingly, eq. (17) is established.

Combining the two sums, we observe

$$\sum_{\frac{n}{2} < k \leq n} (k-1) \mathbb{P}(\mathbb{Y}_n = n - k)$$

$$\leq n \Big(1 - \sum_{1 \leq k < n^{\frac{2}{5}}} \mathbb{P}(\mathbb{Y}_n = n - k) - \sum_{n^{\frac{2}{5}} \leq k \leq \frac{n}{2}} \mathbb{P}(\mathbb{Y}_n = n - k)\Big)$$

$$= n\Big(1 - \big(1 - \alpha n^{-\frac{1}{5}} + o(n^{-\frac{1}{2}})\big) - \big(\alpha n^{-\frac{1}{5}} + o(n^{-\frac{1}{2}})\big)\Big)$$

$$= o(n^{\frac{1}{2}}), \qquad \text{as } n \to \infty.$$

∎

Now we are in position to compute

$$\mathbb{E}[\mathbb{Y}_n] = \sum_{k=1}^{n}(n-k)\mathbb{P}(\mathbb{Y}_n = n-k)$$

$$= n - 1 - \sum_{k=1}^{n}(k-1)\mathbb{P}(\mathbb{Y}_n = n-k)$$

$$= n - 1 - \sum_{1 \le k \le \frac{n}{2}}(k-1)\mathbb{P}(\mathbb{Y}_n = n-k) - \sum_{\frac{n}{2} < k \le n}(k-1)\mathbb{P}(\mathbb{Y}_n = n-k)$$

$$= n - 1 - \alpha n^{\frac{1}{2}}\bigl(1 + o(1)\bigr) - o(n^{\frac{1}{2}})$$

$$= n - \alpha n^{\frac{1}{2}}\bigl(1 + o(1)\bigr).$$

As for the variance,

$$\mathbb{V}[\mathbb{Y}_n] = \mathbb{E}[\mathbb{Y}_n^2] - \mathbb{E}[\mathbb{Y}_n]^2$$

$$= \sum_{k=1}^{n}(n-k)^2 \mathbb{P}(\mathbb{Y}_n = n-k) - \mathbb{E}[\mathbb{Y}_n]^2$$

$$= \sum_{k=1}^{n}(k-1)^2 \mathbb{P}(\mathbb{Y}_n = n-k) + (n-1)^2$$

$$- 2(n-1)\sum_{k=1}^{n}(k-1)\mathbb{P}(\mathbb{Y}_n = n-k) - \mathbb{E}[\mathbb{Y}_n]^2.$$

It is clear from the above computations that $\sum_{k=1}^{n}(k-1)\mathbb{P}(\mathbb{Y}_n = n-k) = \alpha n^{\frac{1}{2}}(1 + o(1))$ and $\mathbb{E}[\mathbb{Y}_n]^2 = (n - \alpha n^{\frac{1}{2}} + o(n^{\frac{1}{2}}))^2$. In this case, we have an analogue of eq. (8)

$$\sum_{n^{\frac{1}{8}} \le k \le \frac{n}{2}} (k-1)^2 \,\mathbb{P}(\mathbb{Y}_n = n-k)$$

$$= \frac{\delta \rho^{\frac{1}{2}} \Phi''(\tau)}{-\Gamma(-\frac{1}{2})\Phi'(\tau)} n^{\frac{3}{2}} \int_0^{\frac{1}{2}} (1-x)^{-\frac{3}{2}} x^{\frac{1}{2}} \mathrm{d}x \,\bigl(1 + o(1)\bigr)\bigl(1 + O(n^{-1})\bigr)$$

$$= \beta n^{\frac{3}{2}}\bigl(1 + o(1)\bigr), \qquad \text{as } n \to \infty,$$

where $\beta = \frac{\delta \rho^{\frac{1}{2}} \Phi''(\tau)}{-\Gamma(-\frac{1}{2})\Phi'(\tau)} \int_0^{\frac{1}{2}}(1-x)^{-\frac{3}{2}} x^{\frac{1}{2}} \mathrm{d}x = (1 - \frac{\pi}{4})\alpha$. Following the same line of arguments, we obtain $\sum_{k=1}^{n}(k-1)^2 \mathbb{P}(\mathbb{Y}_n = n-k) = \beta n^{\frac{3}{2}}(1 + o(1))$. As a result

$$\mathbb{V}[\mathbb{Y}_n] = \beta n^{\frac{3}{2}}(1 + o(1)) + O(n)$$

$$= \beta n^{\frac{3}{2}}(1 + o(1)),$$

completing the proof of the lemma. □

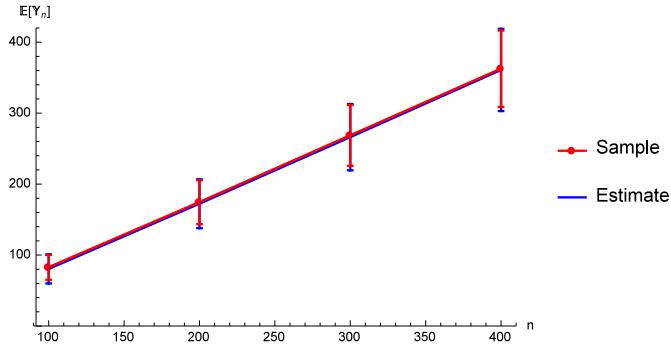

**Fig. 4** The longest rainbow in random RNA secondary structures: we compare the expectation value and standard deviation (blue) with the average length (red) observed in uniformly sampled structures. Minimum arc- and stack-length are $r = \lambda = 1$, $n = 100 \times i$ where $1 \leq i \leq 4$ and sample size is $10^4$, respectively. Error bars represent one standard deviation.

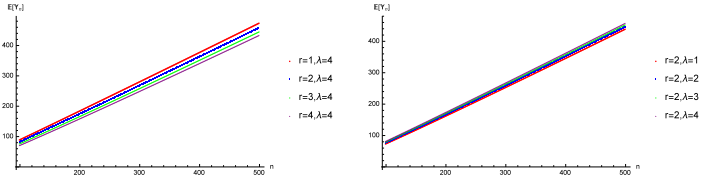

**Fig. 5** The longest rainbow in a random RNA secondary structures: dependency on minimum stack- (LHS) and arc-length (RHS), for $100 \leq n \leq 500$.

In Fig. 4, we contrast our asymptotic estimate of the expectation and the average length of the longest rainbow in random RNA secondary structures.

Fig. 5 shows that the parameter $\alpha$ of the expectation value of the longest rainbow increases, if the minimum stack-length increases or the minimum arc-length decreases.

**Remark:** Lemma 1 shows that the length of the longest rainbow is $n - O(n^{\frac{1}{2}})$ with a standard deviation of $O(n^{\frac{3}{4}})$. As a result, the distribution of $\mathbb{Y}_n$ becomes for larger and larger $n$ more and more concentrated.

**Theorem 3** *We have for any $t > \frac{3}{4}$,*

$$\lim_{n \to \infty} \mathbb{P}(n - \mathbb{Y}_n \geq \Omega(n^t)) = 0 \qquad (18)$$

*and for any $k = o(n)$*

$$\lim_{n \to \infty} \mathbb{P}(n - \mathbb{Y}_n = k) = c\, b_k\, \rho^{k-1}, \qquad (19)$$

*where $c = \Phi'(\mathbf{F}(\rho))^{-1}$, $b_k = [x^{k-1}]\Phi'(\mathbf{F}(x))$ and $\Phi(x) = \frac{1}{1-x}$. Consequently the distribution of $n - \mathbb{Y}_n$ a.a.s. converges to a discrete limit law.*

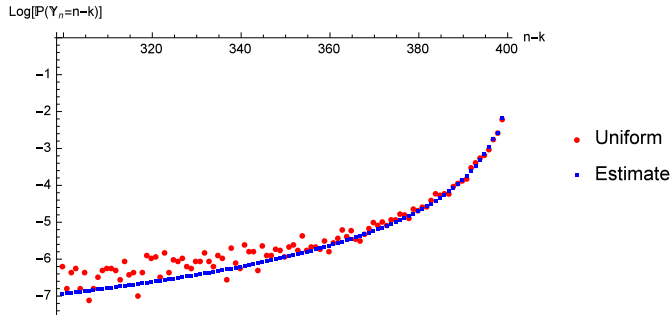

**Fig. 6** The longest rainbow: we contrast the limit distribution (squares) with the distribution (dots) of uniformly sampled structures. Minimum arc- and stack-length are $r = \lambda = 1$, $n = 400$ for a sample size of $10^4$ structures.

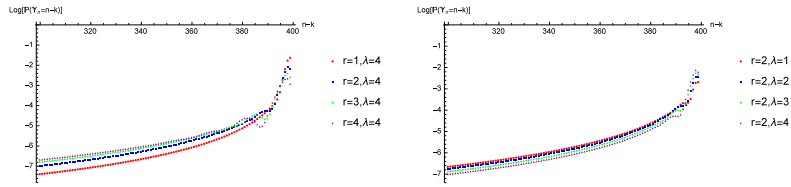

**Fig. 7** The longest rainbow: dependency on minimum stack- (LHS) and arc-length (RHS), where $n = 400$.

*Proof* According to Lemma 1 we have that $\mathbb{Y}_n$ is concentrated at $n - \alpha\, n^{\frac{1}{2}}$ and the variance is $O(n^{\frac{3}{2}})$. Chebyshev's inequality, guarantees

$$\mathbb{P}\Big(\mathbb{E}[\mathbb{Y}_n] - \mathbb{Y}_n \geq a\Big) \leq \frac{\mathbb{V}[\mathbb{Y}_n]}{a^2}.$$

Accordingly, for $a = \Omega(n^t)$ with $t > \frac{3}{4}$, the right hand-side tends to zero as $n$ tends to infinity, whence eq. (18). To establish eq. (19) we inspect that the proof of eq. (11) in Lemma 1 holds for $k = o(n)$ and eq. (19) follows. Eq. (18) implies that a.a.s. we may assume $k = o(n)$ in which case eq. (19) guarantees that $n - \mathbb{Y}_n = k$ satisfies a discrete limit law. □

In Fig. 6, we compare our theoretical result with the distribution of the length of the longest rainbow in uniformly generated structures.

Fig. 7 shows that the decrease of $\mathbb{P}(\mathbb{Y}_n = n - k)$, for increasing $k$, depends on minimum stack- and arc-length.

## 4 The spectrum of rainbows

In the previous section we established that there exists a.s. a unique longest rainbow in an RNA secondary structure. We shall call this rainbow the long rainbow and refer to any other rainbow as short. In this section we study the length-distribution of short rainbows.

To this end we first prove that with high probability we can assume, that any short rainbow is actually finite. That is we show

**Corollary 1** *Given any $\epsilon > 0$, there exists an integer $t(\epsilon)$ such that*

$$\lim_{n \to \infty} \mathbb{P}(\mathbb{Y}_n \geq n - t(\epsilon)) \geq 1 - \epsilon.$$

*In particular, for $r = 1$, $\lambda = 1$, we have*

$$\lim_{n \to \infty} \mathbb{P}(\mathbb{Y}_n \geq n - 100) = 0.808013 \text{ and } \lim_{n \to \infty} \mathbb{P}(\mathbb{Y}_n \geq n - 500) = 0.912911$$

*and in case of $r = 2$, $\lambda = 4$, we have*

$$\lim_{n \to \infty} \mathbb{P}(\mathbb{Y}_n \geq n - 100) = 0.811441 \text{ and } \lim_{n \to \infty} \mathbb{P}(\mathbb{Y}_n \geq n - 500) = 0.913361.$$

*Proof* We observe that $\Phi'(\mathbf{F}(x))$ converges at $\rho$. I.e. for any $\epsilon > 0$, there exists an integer $t(\epsilon)$ such that

$$\sum_{k > t(\epsilon)} c\, b_k \rho^{k-1} < \epsilon,$$

where $c = (\Phi'(\mathbf{F}(x))|_{x=\rho})^{-1}$ and $b_k = [x^{k-1}]\Phi'(\mathbf{F}(x))$. According to Theorem 3, we have

$$\lim_{n \to \infty} \mathbb{P}(\mathbb{Y}_n \geq n - t(\epsilon)) = \lim_{n \to \infty} \sum_{k \leq t(\epsilon)} \mathbb{P}(\mathbb{Y}_n = n - k)$$

$$= \sum_{k \leq t(\epsilon)} c\, b_k\, \rho^{k-1}$$

$$= \sum_{k \geq 1} c\, b_k\, \rho^{k-1} - \sum_{k > t(\epsilon)} c\, b_k \rho^{k-1}$$

$$\geq c\, \Phi'(\mathbf{F}(x))|_{x=\rho} - \epsilon$$

$$= 1 - \epsilon.$$

In case of $r = 1$ and $\lambda = 1$ we obtain $\rho = \frac{1}{3}$ and $c = \frac{1}{9}$ and in case of $r = 2$ and $\lambda = 4$ we have $\rho = 0.540857$ and $c = 0.107902$. This follows by direct computation using Theorem 3. □

In the following, we shall study the distribution of rainbows of finite length. For fixed $k$, let $s_k(n, b)$ denote the number of $r$-canonical secondary structures with minimum arc-length $\lambda$, filtered by the number $b$ of rainbows of length $k$. Let $\mathbf{S}_k(x, u) = \sum_{n,b} s_k(n, b) x^n u^b$ denote the corresponding bivariate generating function.

**Lemma 2** *The bivariate generating function of the number of $r$-canonical secondary structures with minimum arc-length $\lambda$, filtered by rainbows of length $k$, is given by*

$$\mathbf{S}_k(x, u) = \frac{1}{1 - \mathbf{F}(x) - (u - 1)f(k+1)x^{k+1}}.$$

**Remark:** The idea here is to enhance the combinatorial construction underlying the proof of Theorem 1, by marking each rainbow of length $k$. That is, we label each irreducible structure of length $k+1$ using the term $(u-1)f(k+1)x^{k+1}$.

Next we analyze $\mathbb{X}_{k,n}$, the r.v. counting the number of rainbows of length $k$ in a random RNA secondary structure over $n$ nucleotides. By construction, we have

$$\mathbb{P}(\mathbb{X}_{k,n} = b) = \frac{s_k(n,b)}{s(n)} = \frac{[x^n u^b]\mathbf{S}_k(x,u)}{[x^n]\mathbf{S}(x)}.$$

**Theorem 4** *For fixed $k$, $\mathbb{X}_{k,n}$ satisfies the discrete limit law*

$$\lim_{n\to\infty} \mathbb{P}(\mathbb{X}_{k,n} = b) = (b+1)t^b(1-t)^2,$$

*where $\tau = \mathbf{F}(\rho)$ and $t = \frac{f(k+1)\rho^{k+1}}{1-\tau+f(k+1)\rho^{k+1}}$. That is, the limit law of $\mathbb{X}_{k,n}$ is a negative binomial distribution $NB(2,t)$ and the probability generating function of the limit distribution is given by*

$$p_k(u) = \left(\frac{1-t}{1-tu}\right)^2.$$

*Proof* Since $\Phi(x) = \frac{1}{1-x}$ and $h(x,u) := \mathbf{F}(x) + (u-1)f(k+1)x^{k+1}$ have nonnegative coefficients and $h(0,0) = 0$, the composition $\Phi(h(x,u))$ is a well-defined formal power series. In view of Lemma 2, $\mathbf{S}_k(x,u)$ can be expressed as $\mathbf{S}_k(x,u) = \Phi(h(x,u))$.

We verify that $\mathbf{S}_k(x,u)$ has the same dominant singularity $\rho$ as $\mathbf{F}(x)$, by checking that there exists a neighborhood $U$ of 1 such that $h(\rho, u) < 1$ for all $u$ in $U$. As a result, the composition $\mathbf{S}_k(x,u) = \Phi(h(x,u))$ belongs to the subcritical case of singularity analysis (Flajolet and Sedgewick, 2009). Based on the singular expansion of $\mathbf{F}(x)$, we can compute the singular expansion of $h(x,u)$ at $\rho$

$$h(x,u) = \tau + (u-1)f(k+1)\rho^{k+1} + \delta(\rho-x)^{\frac{1}{2}}(1+o(1)),$$

where $\tau = \mathbf{F}(\rho)$. Combining this with the regular expansion of $\Phi(x)$ at $\tau_1 = \tau + (u-1)f(k+1)\rho^{k+1}$

$$\Phi(x) = \Phi(\tau_1) + \Phi'(\tau_1)(x-\tau_1)(1+o(1)),$$

we derive the singular expansion of $\mathbf{S}_k(x,u)$ at $\rho$

$$\mathbf{S}_k(x,u) = \Phi(\tau_1) + \Phi'(\tau_1)\delta(\rho-x)^{\frac{1}{2}}(1+o(1)).$$

The transfer theorem (Flajolet and Sedgewick, 2009), then guarantees

$$[x^n]\mathbf{S}_k(x,u) = \Phi'(\tau_1)\delta\, c_k\, n^{-\frac{3}{2}}\rho^n(1+o(1)).$$

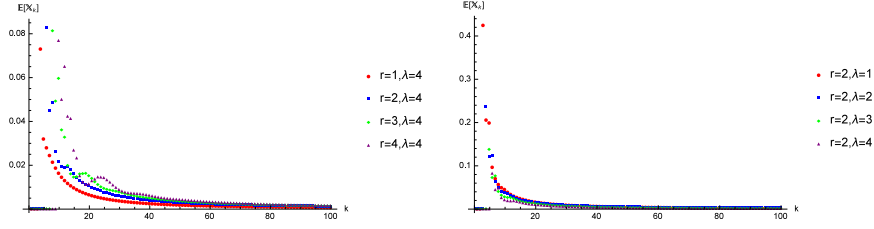

**Fig. 8** The expectation value of rainbows of length $k$ for different minimum stack- and arc-lengths.

Now we are at position to compute

$$p_k(u) = \lim_{n\to\infty} \sum_b \mathbb{P}(\mathbb{X}_{k,n} = b)u^b$$
$$= \lim_{n\to\infty} \frac{\sum_b [x^n u^b]\mathbf{S}_k(x,u)u^b}{[x^n]\mathbf{S}(x)}$$
$$= \lim_{n\to\infty} \frac{[x^n]\mathbf{S}_k(x,u)}{[x^n]\mathbf{S}(x)}$$
$$= \lim_{n\to\infty} \frac{\Phi'(\tau_1)}{\Phi'(\tau)}$$
$$= \left(\frac{1-t}{1-tu}\right)^2,$$

where $t = \frac{f(k+1)\rho^{k+1}}{1-\tau+f(k+1)\rho^{k+1}}$. Exacting the coefficient of $u^b$ in $p_k(u)$, we arrive at

$$\lim_{n\to\infty} \mathbb{P}(\mathbb{X}_{k,n} = b) = [u^b]p_k(u) = \binom{b+1}{b}t^b(1-t)^2.$$

□

**Corollary 2** *For fixed $k$, expectation and variance of $\mathbb{X}_{k,n}$ are asymptotically given by*

$$\lim_{n\to\infty} \mathbb{E}(\mathbb{X}_{k,n}) = \frac{2}{1-\tau}f(k+1)\rho^{k+1}, \quad \lim_{n\to\infty} \mathbb{V}(\mathbb{X}_{k,n}) = \frac{2f(k+1)\rho^{k+1}\left(1-\tau+f(k+1)\rho^{k+1}\right)}{(1-\tau)^2}, \quad (20)$$

*where $\tau = \mathbf{F}(\rho)$.*

Fig. 8 illustrates the dependency of the expectation value of rainbows of length $k$ on minimum stack- and minimum arc-length.

## 5 Discussion

We have shown that the length-spectrum of rainbows in random RNA secondary structures has a gap. By Lemma 1 the longest rainbow is a.a.s. of size $n - O(n^{1/2})$ and Corollary 1 shows that with high probability the second largest rainbow has finite size. In any case, there exists a.a.s. a unique longest rainbow. In Theorem 3 we analyze the limit distribution of the size of the unique longest rainbow and show that it satisfies a discrete limit law. In Theorem 4 we identify the distribution of rainbows of finite size $k$, in the limit of long sequences as a negative binomial.

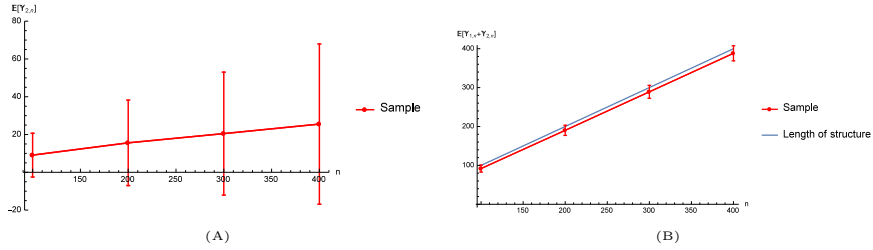

**Fig. 9** We uniformly sample $10^4$ structures with minimum arc- and stack-length $r = \lambda = 1$ and sequence length $n = 100 \times i$ where $1 \leq i \leq 4$. (A) displays the expectation and standard deviation of the length of the second longest rainbow. (B) we display, for the same sample, $k_1 + k_2$, where $k_i$ is the average length of the $i$-th longest rainbow. Error bars represent one standard deviation.

The analysis in Section 3 can be generalized to the lengths of the second and third longest rainbows in uniformly generated structures. One can show that $\mathbb{E}[\mathbb{Y}_{2,n}] = \alpha\, n^{\frac{1}{2}}(1 + o(1))$ and $\mathbb{E}[\mathbb{Y}_{3,n}] = o(n^{\frac{1}{2}})$, where $\mathbb{Y}_{2,n}$ and $\mathbb{Y}_{3,n}$ denote the length of the second and third longest longest rainbow in RNA secondary structures. Suppose the longest rainbow, $\mathbb{Y}_n$, has length $n - k$. Taking out the enclosed irreducible structure, the remaining structure has length $k$. While the rainbow may cut the structure into two distinct intervals of equal orders, the resulting number of structures is far less than the number of structures over a single interval of size $k' = k - o(k)$. In this case, Lemma 1 guarantees that the second longest rainbow has average length $k' + O(k'^{\frac{1}{2}})$ since it is then effectively the longest rainbow of the remaining structure. Therefore, $\mathbb{E}[\mathbb{Y}_{2,n}]$ is $\sum_k (k' + O(k'^{\frac{1}{2}})) \mathbb{P}(\mathbb{Y}_n = n - k) = \sum_k (k + o(k)) \mathbb{P}(\mathbb{Y}_n = n - k)$, which is $\alpha\, n^{\frac{1}{2}}(1 + o(1))$ employing Claim 2 and Claim 3 of Lemma 1. Fig. 9 (A) confirms that the length of the second longest rainbow is $O(n^{\frac{1}{2}})$. Corollary 1 implies that $\mathbb{Y}_{2,n}$ is finite with high probability as $n$ tends to infinity. However we also have $\mathbb{E}[\mathbb{Y}_{2,n}] = O(n^{\frac{1}{2}})$, which means that on a set of measure tending to zero, $\mathbb{Y}_{2,n}$ is infinite. To illustrate this, consider $\mathbb{X}$ with $\mathbb{P}(\mathbb{X} = k) = C_k 4^{-k}$ for $k \geq 1$, where $C_k = \frac{1}{k+1}\binom{2k}{k}$. Then for any $\epsilon > 0$, there exists $k_0$ such that $\mathbb{P}(\mathbb{X} \leq k_0) = 1 - \epsilon$, in other words, $\mathbb{X}$ is finite with high probability. However we have $\mathbb{E}[\mathbb{X}] = \sum_k k \mathbb{P}(\mathbb{X} = k) = \sum_k k C_k 4^{-k} = \infty$.

Our results are connected with the distribution of the $5'$-$3'$ distance in RNA structures, whose finite expectation has been reported in Yoffe et al (2011). Han and Reidys (2012) computed the distribution of these distances proving that they satisfy a discrete limit law. While there is still a set of limit measure zero, composed by structures in which the longest rainbow is not of length $n - k$, for finite $k$, our results show that with high probability we obtain a unique longest rainbow and several short rainbows of finite length. Accordingly, we can provide further insight into the discrete limit law, as the latter does not specify arc-lengths. Fig. 9 (B) shows that the longest and second longest rainbows leave only $o(n^{\frac{1}{2}})$ of nucleotides uncovered. This finding is in accordance with the result of Han and Reidys (2012), who established that the 5'-3' distance is finite.

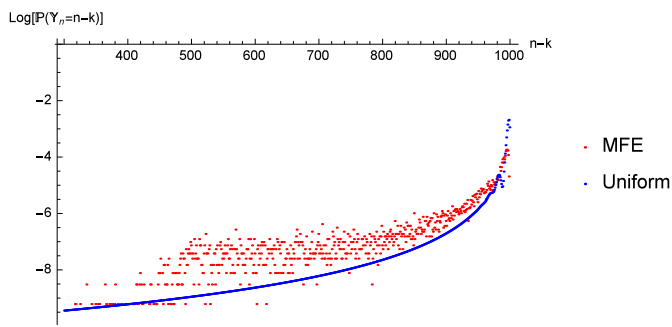

**Fig. 10** The longest rainbow in RNA secondary structures: we compare the limit distribution of random RNA secondary structures having minimum arc- and stack-length four (blue) with the distribution of mfe-structures (red) obtained by ViennaRNA (Lorenz et al, 2011) of $10^4$ random sequences of length 1000.

According to Wexler et al (2007), sparsification achieves linear speed-up if the polymer-zeta property holds. Our results show that in random RNA structures, with high probability, the longest rainbow has almost the length of the sequence. Thus the polymer-zeta property does not hold for RNA secondary structure, unless one considers particular classes of natural RNA structures such as mRNA (Wexler et al, 2007). Having a closer look at the number of stems in random RNA structures, we observe a central limit theorem (thus having an expectation value of order $O(n)$). This suggests that the expected size of stems is $O(1)$ and thus we find $O(1)$ arcs of length $n - O(n^{1/2})$.

Let us put our results into context with mfe-structures. To this end we compare in Fig. 10 the limit distribution of the length of the longest rainbow with that of mfe-structures. We can report that the longest rainbow in mfe-structures satisfies a similar distribution. Closer inspection reveals, that compared to random structures, mfe-structures exhibit fewer rainbows of length between 980 and 1000 and more rainbows of length between 400 and 980. Increasing the minimum stack-size in random structures has the effect that the distribution of lengths of the longest rainbow in random and mfe-structures becomes more and more similar, see Fig. 7 (LHS). This makes sense as mfe-structures are typically form stacks of larger size. As for the expectation value, Fig. 11 shows that the longest rainbow of mfe-structures is also close to, but smaller than that in random structures.

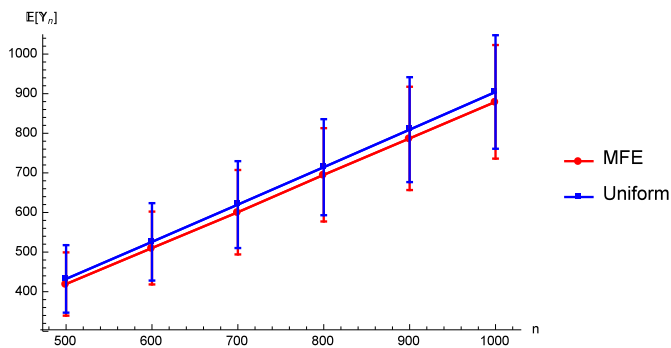

**Fig. 11** Expectation value and standard deviation: we compare the theoretical estimate for random structures having minimum arc- and stack-length four (blue) with the distribution of mfe-structures (red), obtained by ViennaRNA (Lorenz et al, 2011) of $10^4$ random sequences of length $500 \le n \le 1000$.

**Table 1** The probability of having a long rainbow in RNA structures: we contrast our theoretical result in Corollary 1 for $r = 4$ and $\lambda = 4$ and the probabilities obtained from $10^4$ random mfe-structures of length 1000.

| $\mathbb{P}(\mathbb{Y}_{1000} \ge 1000 - k)$ | $k = 100$ | $k = 200$ | $k = 300$ | $k = 400$ | $k = 500$ |
|---|---|---|---|---|---|
| mfe | 0.6333 | 0.7779 | 0.8574 | 0.9207 | 0.9775 |
| uniform | 0.7179 | 0.7936 | 0.8295 | 0.8514 | 0.8666 |

We have shown in Corollary 1 that random structures exhibit with high probability a longest rainbow of size $n - k$, for finite $k$. In Table 1 we study this phenomenon in mfe-structures. In fact we observe that this probability is higher in mfe-structures than in random structures, indicating that for mfe-structures the gap in the sequence of length of rainbows is more pronounced.

In Fig. 12 we display that eq. (20) provides a good approximation for the expected number of short rainbows of length $\ge 25$ in mfe-structures. However, mfe-structures have fewer rainbows of length between 5 and 15 and more rainbows of length $\ge 15$ than random structures. As we observed in the context of the length distribution of the longest rainbow, increasing minimum stack-size in random structures results in a better and better approximation of short rainbows in mfe-structures. This seems plausible, as mfe-structures, in order to achieve minimum energy, tend to form long stems.

Finally, we discuss our findings in the context of rainbows observed in structures contained in RNA databases (from the RCSB PDB database (Berman et al, 2000) and the comparative RNA web (CRW) site (Cannone et al, 2002)). The observed average ratio of the length of the longest rainbow relative to the length of the sequence varies with different RNA families. For tRNA (76–90 nt) this ratio is $0.928(\pm 0.048)$, implying a long rainbow. This is a result of the fact that tRNA typically forms the cloverleaf structure (Kim et al, 1974;

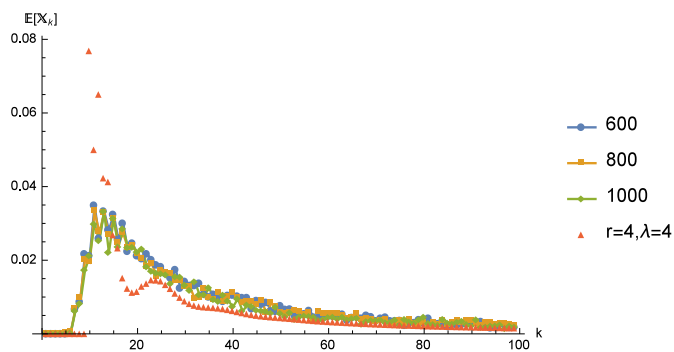

**Fig. 12** The expected number of rainbows of length $k$ in RNA secondary structures: we contrast eq. (20) for random structures having minimum arc- and stack-length four (orange) and mfe-structures, of $10^4$ random sequences of lengths $n = 600, 800, 1000$, respectively.

Robertus et al, 1974). A similar ratio is also observed for transfer-messenger RNA (tmRNA, 300–400 nt), 5S rRNA (120 nt). 23S rRNAs of *Escherichia coli* and *Thermus thermophilus* (2904 nt) exhibit a ratio of 0.999. However 16S rRNAs of the same species (1542 nt) has a ratio of 0.584 for the longest rainbow and a ratio 0.308 for the second longest (Woese et al, 1980). This shows a general tendency of natural RNA structures to have a unique longest rainbow but there are exceptions: specific functionalities lead to structures having a small number of long rainbows.

As for future work we are concerned with the implications of the results of this paper for the entire arc-spectrum of RNA secondary structures. We argue here that, for $n$ sufficiently large, using the fact that arcs are non-crossing, we can employ the results on the longest rainbow in order to compute the entire arc-spectrum in a recursive manner. Namely, once the stack concerning the longest rainbow is removed, we obtain an induced, nested, reducible RNA secondary structure. With respect to this structure we then iterate the argument working our way from top to bottom.

## 6 Acknowledgments.

We would like to thank the reviewers for their comments and suggestions and specifically for pointing out a gap in the proof of Lemma 1. We gratefully acknowledge the help of Kevin Shinpaugh and the computational support team at BI. Many thanks to Christopher L. Barrett and Henning Mortveit for discussions. The second author is a Thermo Fisher Scientific Fellow in Advanced Systems for Information Biology and acknowledges their support of this work.